\documentclass[12pt]{amsart}
\usepackage{amsmath,amssymb}

\newcommand{\Z}{{\mathbb Z}}

\newcommand{\N}{{\mathbb N}}
\newcommand{\eps}{\varepsilon}
\newtheorem{thm}{Theorem}

\newtheorem{conj}{Conjecture}

\newtheorem{lem}{Lemma}

\author[S. Akiyama]{Shigeki Akiyama}
\address[S. Akiyama]{Institute of Mathematics, University of Tsukuba, 1-1-1 Tennodai, Tsukuba, Ibaraki, 305-8571 Japan}
\email{akiyama@math.tsukuba.ac.jp}

\author[B. Mance]{Bill Mance}
\address[B. Mance]{Uniwersytet im. Adama Mickiewicza w Poznaniu,
  Collegium Mathematicum, ul. Umultowska 87, 61-614 Pozna\'{n}, Poland}
\email{william.mance@amu.edu.pl}

\thanks{The first author was supported by JSPS grants in aid 24K06662. The second author was supported by grant 2019/34/E/ST1/00082 for the project ``Set theoretic methods in dynamics and number theory,'' NCN (The National Science Centre of Poland).}

\begin{document}
\title{Discretized Rotation with fixed initial points}
\date{}
\maketitle

\begin{abstract}
Motivated by comparison of discrete and continuous dynamical systems,
we expect that bijective discretized rotations approximate rotation 
dynamics well. 
We examine this hypothesis by
a simple but notorious model of discretized rotation
for which little is known. 
We prove that if $\max\{|a_0|,|a_1|\}\le 10$ and $\lambda\in\ ]-2,2[$,
then the sequence defined by
$$
0 \le a_{n+2} +\lambda a_{n+1}+a_n<1
$$
is periodic.
\end{abstract}

\section{Introduction}

Motivated by comparison of discretized and continuous systems, 
the discretized rotation has attracted many researchers as a toy model.
Our target mapping is the composition of rotation and some
cut-off operation
to fit into the discrete space. There are many studies on 
this type of mappings and several embeddings of the problem in
different spaces (torus or $p$-adic space) are discussed, 
see \cite{Lowenstein-Vivaldi:98, Lowenstein-Hatjispyros-Vivaldi,
BLPV, Bosio-Vivaldi:00, Vivaldi:06}.

To expect properties similar to a rotation, our
discretization should be a bijective map, see \cite
{KKPV, Berthe-Nouvel, Roussillon-Coeurjolly, Hannusch-Pethoe:23}.
Since the orbit of the rotation is bounded, we may
expect that
the orbit is periodic for the bijective discretized rotation. 
In this paper we study a simple case:

\begin{conj}
\label{SRS}
For any  $\lambda\in\ ]-2,2[$
and any
initial value $(a_0,a_1)\in \Z^2$, the integer sequence $(a_n)_{n=0,1,\dots}$ defined by
\begin{equation}
\label{SRS0}
0 \le a_{n+2} +\lambda a_{n+1}+a_n<1
\end{equation}
is periodic.
\end{conj}

By symmetry of the definition, 
the associated mapping is bijective and was studied by J. H. Lowenstein, S. Hatjispyros, and F. Vivaldi in \cite{Lowenstein-Hatjispyros-Vivaldi} and 
S. Akiyama, H. Brunotte, A. Peth\"{o}, and J. M. Thuswaldner in \cite{Akiyama-Brunotte-Pethoe-Thuswaldner:05}. 
The map $(a_0,a_1) \to (a_1,a_2)$ is written as
$$
\begin{pmatrix} x \\ y \end{pmatrix}
\mapsto \begin{pmatrix} y \\ -\lfloor x+\lambda y\rfloor\end{pmatrix}.
$$
This can be rewritten as
$$
\begin{pmatrix} x \\ y \end{pmatrix}
\mapsto \begin{pmatrix} 0 & 1\\
                    -1& -\lambda \end{pmatrix}
\begin{pmatrix} x \\ y \end{pmatrix} + 
\begin{pmatrix} 0 \\ \eps \end{pmatrix}
$$
with $\eps\in [0,1[$.  The associated matrix has two
characteristic roots of modulus $1$, complex conjugate to each other.
Therefore if we fix $\lambda\in\ ]-2,2[$, then 
this map is a bijective discretized rotation acting on a lattice. 
It is clear that 
the orbit is periodic if and only if it is bounded. 
It is also equivalent to pure periodicity because the map is a bijection.
Therefore Conjecture \ref{SRS} asks if 
the accumulation of error vectors $(0,\eps)^T$ cancels the effect of rotation and we obtain
a bounded sequence. We believe that this heuristic argument 
supports Conjecture~\ref{SRS}, but little is known around this Conjecture. 

Conjecture~\ref{SRS} is known to be
 valid for 11 values of $\lambda$
$$
\lambda=0,\ \pm 1,\ \frac{\pm 1 \pm \sqrt{5}}{2},\ \pm \sqrt{2},\ \pm \sqrt{3}.
$$
The cases $0,\pm 1$ are trivial. However, the proofs of the other 
cases are not easy. 
J. H. Lowenstein, S. Hatjispyros, and F. Vivaldi \cite{Lowenstein-Hatjispyros-Vivaldi} transferred 
the problem to the map on the 
torus and gave a computer assisted proof for $\lambda=(1-\sqrt{5})/2$.
A number theoretic proof for $\lambda=(1+\sqrt{5})/2$
was given in by S. Akiyama, H. Brunotte, A. Peth\"{o}, and W. Steiner \cite{Akiyama-Brunotte-Pethoe-Steiner:06}.
In \cite{Akiyama-Brunotte-Pethoe-Steiner:07} they gave a proof
by using the self-inducing structure of the associated domain exchanges,
which makes the proofs simpler and they didn't require computer assistance.
The embedding to the torus makes an interesting bridge to another problem
on piecewise isometry, see H. Bruin, A. Lambert, G. Poggiaspalla, and S. Vaienti \cite{BLPV} and  K. Kouptsov, J. H. Lowenstein, and F. Vivaldi \cite{Kouptsov-Lowenstein-Vivaldi:02}.
An exposition by S. Akiyama and E. Harriss can be found in \cite{Akiyama-Harriss:13}.

There has been little progress on this problem for a fixed $\lambda$ 
since then. 
S. Akiyama and A. Peth\"{o} \cite{Akiyama-Petho:12} proved that 
there exists infinitely many periodic orbits in Conjecture \ref{SRS}. 
This is a weak statement, but we believe it is far from trivial.

In this paper, we propose a totally different approach. We
shall fix an initial vector $(a_0,a_1)\in \Z^2$
and study which $\lambda \in\ ]-2,2[$ satisfy Conjecture~\ref{SRS}. 
This approach is wild as we are discussing totally 
different dynamical systems of different rotation angles 
at the same time, but interestingly we can obtain several new observations. 

Let us assume that Conjecture \ref{SRS} 
and its infinitesimal small perturbations 
are valid.
More precisely, 
we assume not only Conjecture \ref{SRS} but also
the periodicity of the limit of the integer sequence
$(b_n(\eps))_{n=0,1,\dots}$ defined by
\begin{equation}
\label{SRSA1}
0 \le b_{n+2}(\eps) +\left(\frac{p}{q}+ \eps\right) b_{n+1}(\eps)+b_n(\eps)
<1
\end{equation}
with $(b_0(\eps),b_1(\eps))=(a_0,a_1)$, as $\eps \to +0$ by 
the product topology of the sequence space $\Z^{\N}$
for any irreducible fraction $p/q\in ]-2,2[$.
Since $b_n(\eps)\in \Z$, the limit is
well-defined. Further we also
assume the periodicity of the limit for the integer sequence
$(c_n(\eps))_{n=0,1,\dots}$ defined by
\begin{equation}
\label{SRSA3}
0 \le c_{n+2}(\eps) +\left(\frac{p}{q}- \eps\right) c_{n+1}(\eps)+c_n(\eps)
<1
\end{equation}
with $(c_0(\eps),c_1(\eps))=(a_0,a_1)$.
To be concrete, since the inequality (\ref{SRSA1}) with $\eps \to +0$
gives rise to the map
$$
\begin{pmatrix} x \\ y \end{pmatrix}
\mapsto 
\begin{cases}
\begin{pmatrix}
y \\ 1+ \left\lceil - \frac{p}{q}  y -x \right\rceil
\end{pmatrix}
&
y <0 \land y\equiv 0 \pmod{q} \\
\begin{pmatrix} 
y \\ \left\lceil - \frac{p}{q}  y -x \right\rceil
\end{pmatrix}
& \text{ otherwise},
\end{cases}
$$
we assume that orbit of this map starting from any $(a_0,a_1)\in \Z^2$
is periodic. Similarly for the inequality (\ref{SRSA3}), 
we assume the periodicity of the orbits of
$$
\begin{pmatrix} x \\ y \end{pmatrix}
\mapsto
\begin{cases}
\begin{pmatrix} 
y \\ 1+ \left\lceil - \frac{p}{q}  y -x \right\rceil
\end{pmatrix}
&
y >0 \land y\equiv 0 \pmod{q} \\
\begin{pmatrix} 
y \\ \left\lceil - \frac{p}{q}  y -x \right\rceil
\end{pmatrix}
& \text{ otherwise}.
\end{cases}
$$
For simplicity, we write
$$
0 \le a_{n+2} +\left(\frac{p}{q}\pm 0\right) a_{n+1}+a_n<1
$$
for (\ref{SRSA1}) or (\ref{SRSA3}). 
This notation appears in section
\ref{Lems} as well.
Assumptions (\ref{SRSA1}) and (\ref{SRSA3})
are minor modifications of the original conjecture and are expected to be 
true from the above heuristic argument. 

We try to see how the period of (\ref{SRS0})
changes when $\lambda$ varies within $]-2,2[$. 
%
%
Our observations are as follows.

\begin{thm}
\label{Main}
Assume Conjecture \ref{SRS} and 
the periodicity (\ref{SRSA1}) and
(\ref{SRSA3}). Then 
for any $(a_0,a_1)\in \Z^2$ and $\eps>0$ there exists a finite number of 
non empty intervals $I_j$ such that
\begin{itemize}
\item $[-2+\eps,2[= \bigsqcup_{j} I_j$
\item For any $\lambda\in I_j$, the period of $(a_n)$ by (\ref{SRS0})
is identical.
\item If $\lambda\in I_j$ and $\lambda'\in I_j'$ and $j\neq j'$, the period of 
(\ref{SRS0}) for $\lambda$ and the one for $\lambda'$ are different.
\end{itemize}
\end{thm}

Here the interval $I_j$ is closed, open or semi-open. 
It can be a singleton.

\begin{thm}
\label{Main2}
Assume  Conjecture \ref{SRS} and the periodicity of (\ref{SRSA1}) and
(\ref{SRSA3}).
If $(a_0,a_1)=(m,m)$ with $m\ge 0$, then there exists a
finite number of non empty intervals $I_j$ such that
\begin{itemize}
\item $]-2,2[= \bigsqcup_{j} I_j$
\item For any $\lambda\in I_j$, the period of $(a_n)$ by (\ref{SRS0})
is identical.
\item If $\lambda\in I_j$ and $\lambda'\in I_j'$ and $j\neq j'$, the period of 
(\ref{SRS0})
for $\lambda$ and the one for $\lambda'$ are different.
\end{itemize}
\end{thm}

One can also see that 
if $a_0\neq a_1$ then a
countably infinite collection of intervals having such properties
is required to cover $]-2,2[$. 
Moreover, we can make explicit the shape of the
intervals which cover $]-2,-2+\eps[$ for a small $\eps>0$ in 
Theorem \ref{Main3}.
Theorem \ref{Main} covers the interval $[-2+\eps,2[$ and 
we can make $I_j$ explicit when $a_0,a_1$ are 
fixed by symbolic computation. Thus we can obtain our main theorem: an unconditional result:

\begin{thm}
\label{Small}
If $\max\{|a_0|,|a_1|\}\le 10$ and  $\lambda\in\ ]-2,2[$,
then the sequence defined by
$$
0 \le a_{n+2} +\lambda a_{n+1}+a_n<1
$$
is periodic.
\end{thm}

The upper bound $10$ can likely be improved but it depends on the computational cost. 
For a discussion, see the closing remarks.
We can also obtain another interesting unconditional result:

\begin{thm}
\label{2}
For any $(a_0,a_1)\in \Z^2$ there exists $\eps>0$ so that
if $2-\eps<\lambda<2$, then the sequence defined by
$$
0 \le a_{n+2} +\lambda a_{n+1}+a_n<1
$$
is periodic.
\end{thm}
Furthermore, we have
\begin{thm}
\label{2}
For any $(a_0,a_1)\in \Z^2$ there exists $\eps>0$ so that
if $-2<\lambda<-2+\eps$, then the sequence defined by
$$
0 \le a_{n+2} +\lambda a_{n+1}+a_n<1
$$
is periodic. A value of $\eps$ is explicitly given in Theorem~\ref{Main3}, where a full description of the structure of the cycles on this interval are given.
\end{thm}

\section{Lemmas: around $\pm 2$}
\label{Lems}

We employ the same abbreviation as (\ref{SRSA1}) 
and (\ref{SRSA3}) in the previous section.
We wish to study the integer sequence defined by
\begin{equation}
\label{GDef}
0\le a_{n+2} +(2-0) a_{n+1}+a_n<1.
\end{equation} 
Then we obtain
$$a_{n+2}= \begin{cases} -a_n -2a_{n+1} & a_{n+1}\le 0 \cr 
                     -a_n- 2a_{n+1}+1  & a_{n+1}>0\end{cases}.
$$
This defines a map
$$
G((x,y))= \begin{cases} (y, -x-2y) & y \le 0\cr
                         (y, -x-2y+1) & y>0\end{cases}
$$
from $\Z^2$ to itself, which is bijective by the symmetry of 
(\ref{GDef}). 
Also we are interested in the sequence 
$$
0\le a_{n+2} +(-2+0) a_{n+1}+a_n<1.
$$

Then
\begin{equation}
\label{Hdef}
a_{n+2}= \begin{cases} -a_n+2a_{n+1} & a_{n+1}\ge 0 \cr 
                     -a_n+2a_{n+1}+1  & a_{n+1}<0
\end{cases}
\end{equation}
Therefore we define another bijection
$$
H((x,y))= \begin{cases} (y, -x+2y) & y \ge 0\cr
                         (y, -x+2y+1) & y<0\end{cases}.
$$
The definitions of $G$ and $H$ are similar but their behaviors are very different. Clearly $H((m,m))=(m,m)$ for $m\ge 0$.

\begin{lem}
\label{H}
For $(x,y)\in \Z^2$ with $(x,y) \not \in \{ (m,m)\ | m=0,1,\dots\}$ 
the orbit of $H$ is unbounded. 
\end{lem}

\begin{proof}
Letting $k((x,y))=x-y$, we have
$$
k(H((x,y))) = \begin{cases} k((x,y)) & \text{ when } y \ge 0 \cr
                            k((x,y))-1 & \text{ when } y<0\end{cases}.
$$
This implies $k((x,y))\ge k(H((x,y)))$.

Therefore if $k((x,y))<0$, then
$k(H^n((x,y)))<0$ holds for all $n=0,1,\dots$. Writing 
$(a_n,a_{n+1})= H^n((x,y))$, we have $a_n<a_{n+1}$ for all $n=0,1,\dots$. Thus
$a_n \to \infty$ as $n\to \infty$
and the orbit of $H$
is unbounded. The proof is finished in this case.

Take $(x,y)\not \in \{ (m,m)\ | m=0,1,\dots\}$. Since $H$ is a bijection, 
and $H((m,m))=(m,m)$ for $m\ge 0$, 
there exists no $n\in \N$ with $H^n(x,y)\in \{ (m,m)\ | m=0,1,\dots\}$.
Thus if there exists $n\in \N$ with 
$k(H^n(x,y))=0$, then $a_n=a_{n+1}<0$. 
Thus $k(H^{n+1}((x,y))<k(H^n(x,y))=0$ and the proof is finished 
using $H^{n+1}((x,y))$ 
instead of $(x,y)$. 

Finally we assume 
that $k(H^n((x,y)))>0$ for all $n=0,1,\dots$. 
In this case, 
we have $a_n>a_{n+1}$ for all $n$. Thus $a_n$ diverges to $-\infty$.
\end{proof}

\begin{lem}
\label{G}
For any $(x,y)\in \Z^2$, the orbit of $G$ is periodic.
\end{lem}

\begin{proof}
In this case, the key quantity is $k((x,y))=x+y$. The orbit from
$(m,-m)$ with $0<m\in \Z$ must be symmetric, since $G((m,-m))=(-m,m)$
and this implies $G^{-n}((m,-m))=\tau(G^n(-m,m))$ for $n\in \N$
by the symmetry of (\ref{GDef}). 
Here $\tau((x,y))=(y,x)$. 
The symmetric orbit is periodic if and only if
it is doubly symmetric.
We partition the set
$\Z^2 \setminus \{(m,-m)\ |\  m\in \Z\}$
into four parts. 
\begin{align*}
A&=\{(x,y)\in \Z^2\ |\ x+y >0, y> 0\}\\
B&=\{(x,y)\in \Z^2\ |\ x+y <0, y> 0\}\\
C&=\{(x,y)\in \Z^2\ |\ x+y <0, y\le 0\}\\
D&=\{(x,y)\in \Z^2\ |\ x+y >0, y\le 0\}.
\end{align*}
We claim that $G$ acts in the following way:
\begin{align*}
&G^2(A)\subset A \cup \{(-m,m)\ |\ m>0\}\\
&G^2(C)\subset C \cup \{(m,-m)\ |\ m>0\}\\
&G(B)\subset D \cup A\\
&G(D)\subset B \cup C.
\end{align*}
Indeed, 
the map $G^2$ is written as
\begin{equation}
\label{Dec}
G^2((x,y))=(-x-2y+1,2x+3y-2)
\end{equation}
when $(x,y)\in B \cup A$ and the first inclusion follows.
Similarly the second one is derived from
\begin{equation}
\label{Inc}
G^2((x,y))=(-x-2y,2x+3y+1)
\end{equation}
when $(x,y)\in C \cup D$. Other two inclusions are direct consequences 
of the definition of $G$. The claim is proved.

Further, we see that the orbits of 
points in $B$ and $D$ eventually fall into $A$ or $C$. 
 If there exists an infinite orbit switching 
between $B$ and $D$, then we find
infinite iteration of the map $B\to D\to B$ of the shape (\ref{Dec}). 
We obtain
$$
G^{2n}((x,y))=(n^2-(2 n - 1) x - 2 n y, -n (n + 1) + 2 n x + (2 n + 1) y).
$$
The first term 
$n^2-(2 n - 1) x - 2 n y\to \infty$ as $n\to \infty$. 
Since the first coordinate of $B$ is not positive, we get a contradiction. 
This shows that the number of iterates of $B\to D\to B$ is bounded.

For $(x,y)\in A$ we have $k((x,y))=k(G^2((x,y)))-1$ from (\ref{Dec})
and the orbit of $G^2$ must reach a point of the form $(-m,m)$ with $m>0$.
Similarly for $(x,y)\in C$ we have $k((x,y))=k(G^2((x,y)))+1$ from (\ref{Inc})
and the orbit of $G^2$ must reach a point of the form $(m,-m)$ with $m>0$.
Since $G((m,-m))=(-m,m)$ for $m>0$ gives a center of
symmetry of the orbit of $G$
and the orbit starting from $(-m,m)$ must follow the same dynamics. 
Thus we have shown that
every orbit is symmetric, and its forward orbit must fall into another 
center of symmetry
(possibly the same one). 
This shows that the orbit must be doubly symmetric and therefore
periodic. 

The period is classified into two types:
even cycles of the form 
$$(m,-m,m, c_1,\dots, c_t , \ell, -\ell, \ell, c_t, \dots c_1)^{\infty}
$$
with $m\neq \ell$ and odd cycles of the form
$$(m,-m,m, c_1,\dots, c_t, \ell, \ell, c_t, \dots c_1)^{\infty}.
$$
In particular, all cycles are cyclic palindromes.
\end{proof}

\section{Proof of Theorems  \ref{Main}, \ref{Main2} and \ref{2}}

By the periodicity of (\ref{SRS0}), 
a fixed initial value $(a_0,a_1)$ and any $\lambda\in\ ]-2,2[$, 
there exists $n\ge 1$ such that $a_n=a_0$, $a_{n+1}=a_1$, i.e.,
$$
0\le a_{i+2} + \lambda a_{i+1}+a_{i}<1
$$
for $i=0,1,\dots, n-1$ and $(a_n,a_{n+1})=(a_0,a_1)$. 
We consider the system of $2n$ inequalities 
$$
0\le a_{i+2} + x a_{i+1}+a_{i}<1
$$
with a variable $x\in\ ]-2,2[$.
When $a_{i+1}\neq 0$, every inequality gives an interval $[y, 2[$, 
 $]y, 2[$, $]-2, y]$, or $]-2, y[$ with some $y$.
Therefore the intersection of $2n$ inequalities is a non empty interval. Since the coefficients of the inequalities are integers, the resulting interval has rational end points. The interval can be open, closed, or semi-open. Note that
it can be a singleton.
In this way for a fixed $(a_0,a_1)$, 
we obtain 
non-empty rational intervals $I$ 
which correspond 
one to one to a cycle $(b_0,b_1,\dots, b_{n-1})^{\infty}$  with $(b_0,b_1)=(a_0,a_1)$ generated by $\lambda\in I_j$. 

This gives a countable partition :
\begin{equation}
\label{Part}
]-2,2[ = \bigsqcup_{j=1}^{\infty} I_j.
\end{equation}
We are interested in the accumulation points of $E:=\bigcup_{j=1}^{\infty} \partial(I_j)$: the set of end points of $I_j$, in the ambient space $[-2,2]$. 

We claim
that there exists no accumulation point of $E$ 
within $]-2,2[$. As the end points
of $I_j$ are rational, every irrational point in $]-2,2[$ is
an inner point of some $I_j$. Thus an irrational point is not an 
accumulation point of $E$. Take a rational number $p/q\in\ ]-2,2[$.
Again by (\ref{SRSA1}) for $\lambda=p/q$, 
there exists an interval $I_j$ which contains $p/q$ and
its associated period is
$(b_0,b_1,\dots, b_{n-1})^{\infty}$ with $(b_0,b_1)=(a_0,a_1)$.
Let us compare two sequences. The first is defined by
$$
0 \le b_{n+2} +\left(\frac{p}{q}+ 0\right) b_{n+1}+b_n<1
$$
with $(b_0,b_1)=(a_0,a_1)$. The second is
$$
0\le d_{n+2}+\left(\frac pq+v\right) d_{n+1} + d_n <1
$$
for a small $v$, with $(d_0,d_1)=(a_0,a_1)$. 
By (\ref{SRSA1}), $(b_n)$ is periodic and we put $V=\max_{i} |b_i|$.
By induction, using  $b_{n+2} +\frac{p}{q} b_{n+1}+b_n\in (1/q)\Z$, 
we see that the two sequences $(b_i)$ and $(d_i)$
are identical if $0<v< 1/(2qV)$. 
Conjecture \ref{SRS} guarantees that for a fixed initial vector
$(a_0,a_1)\in \Z^2$, we have a countable partition (\ref{Part}).
The above discussion shows that if $0<v<1/(2qV)$, then the orbits are equal to $(b_n)$.
This implies that the point $p/q$ is not a right accumulation point of $E$.
Similarly using (\ref{SRSA3}), $p/q$ is not a left accumulation point of $E$
as well. Our claim is proved.


To prove Theorem \ref{Main}, it suffices to show that $2$ is not an accumulation point of $E$. 
Let us compare two sequences. The first is defined by
$$
0 \le b_{n+2} +\left(2- 0\right) b_{n+1}+b_n<1
$$
with $(b_0,b_1)=(a_0,a_1)$, i.e., $(b_n,b_{n-1})=G^n((a_0,a_1))$
by the function $G$ defined in the previous section.
By Lemma \ref{G}, $(b_n)$ is periodic and we put $V=\max_i |b_i|$.
The second is
$$
0\le c_{n+2}+\left(2-v\right) c_{n+1} + c_n <1
$$
for a small $v$, with $(c_0,c_1)=(a_0,a_1)$. If we take $0<v<1/(2V)$, 
the sequence $(b_n)$ and $(c_n)$ are identical by induction. Therefore
if $\lambda\in\ ]2-v,2[$, then it corresponds to the same cycle $(b_i)$.
This shows that $2$ is not an accumulation point of $E$. 
The proof of Theorem \ref{Main} is finished. 
Since Lemma \ref{G} is an unconditional statement, we immediately 
obtain Theorem \ref{2} as well.

Let us consider $(a_0,a_1)=(m,m)$ with $m>0$. The sequence
defined by
$$
0\le a_{n+2}+ (-2 +v) a_{n+1}+ a_n<1
$$
is the constant cycle $(m)^{\infty}$ if $0\le v<1/m$. Therefore in this case $-2$ is covered by an interval $]-2,-2+1/m[$. Because of Theorem \ref{Main}, there are 
no accumulation points of $E$
in $[-2,2]$ and the proof of Theorem \ref{Main2}
is finished. 

Lemma \ref{H} shows that the point $-2$ is indeed an accumulation point of $E$
for an initial vector $(a_0,a_1)$ with $a_0\neq a_1$. 
This will be made precise in Theorem \ref{Main3}.

\section{Theorem \ref{Main3}: around $\lambda=-2$}

Let $T_n=n(n+1)/2$ be the n'th triangular number.
Given two integers $s,d$ with $0\le s<d$ and a positive integer $k$, let
$$
A=(s,s+d,s+2d,\cdots,s+(k-1)d,s+kd)
$$
Next, we start adding differences of triangular numbers:
$$
B=(s+kd+(T_d-T_{d-1}),s+kd+(T_d-T_{d-2}),\ldots,s+kd+(T_d-T_{1}),s+kd+(T_d-T_{0}))
$$
Finally, consider the descent into negative integers
$$
C=(s-(T_d-T_{d-1}),s-(T_d-T_{d-2}),\ldots,s-(T_d-T_0)).
$$
Define a cycle
$$
\mathcal{C}_{s,d,k}:=(AB\bar{B}\bar{A}C\bar{C})^{\infty}.
$$
The bar reverses the order. Then we have a

\begin{lem}
\label{Hex}
If $k\ge (T_d-s)/d$, then $\mathcal{C}_{s,d,k}$ corresponds to 
the interval
$$
\left[-2+\frac 1{s+(k+1)d}, -2+\frac 1{s+kd} \right[.
$$
\end{lem}

\begin{proof}
We denote $A=(a_i)$, $B=(b_j)$, $C=(c_j)$ with $i=0,1,\dots, k$ and $j=1,\dots, d$. We notice that their definitions can be overlappingly 
extended to the cycle $\mathcal{C}_{s,d,k}$, i.e., 
$a_{k+1}=s+(k+1)d=s+kd+T_d-T_{d-1}=b_1$, $b_0=s+kd+(T_d-T_d)=a_k$, 
$b_{d+1}=s+kd+(T_d-T_{-1})=b_d$, 
$a_{-1}=s-d=s-(T_d-T_{d-1})=c_1$, 
$c_0=s-(T_d-T_d)=a_0$, $c_{d+1}=s-(T_d-T_{-1})=s-T_d=c_d$. 
Therefore we only have to check that
\begin{align*}
0\le& s+(i-1)d +(-2+\eps)(s+id) + s+(i+1)d<1\\
0\le& s+kd +(T_d-T_{d-j+1}) +  (-2+\eps) (s+kd +(T_d-T_{d-j}))\\
&+ s+kd +(T_d-T_{d-j-1})<1\\
0\le& s-(T_d-T_{d-j+1}) +  (-2+\eps) (s-(T_d-T_{d-j}))+ 
s-(T_d-T_{d-j-1})<1
\end{align*}
with $i=0,1,\dots, k$ and $j=1,\dots d$ is valid for $$\eps \in [1/(s+(k+1)d),1/(s+kd)[.$$ 
Using $T_{j-1}-2T_j+T_{j+1}=1$, they simplify  to
\begin{align*}
0&\le \eps (s+ kd) <1\\
1&\le \eps (s+(k+1)d)<2\\ 
1&\le \eps (s+kd +T_d)<2\\
-1&\le  \eps (s-d) <0\\
-1&\le  \eps (s-T_d) <0.
\end{align*}
and using $T_d\le s+kd$, we see that the truth of all these conditions is equivalent to
$$
\frac 1{s+(k+1)d} \le \eps < \frac 1{s+kd}.
$$
\end{proof}

We wish to assign every point $(a_0,a_1) \in \Z^2$ a label $(s,d)$ 
with $0\le s\le d$
which will partition $\Z^2$ into an equivalence class based on this label:
\begin{enumerate}
\item If $0 \leq a_0 < a_1$, then $d=a_1-a_0$ and $s=a_0-d\lfloor a_0/d\rfloor$. 
\item If $0 \leq a_1 < a_0$, then $d=a_0-a_1$ and $s=a_1-d\lfloor a_1/d\rfloor$. 
\item If $a_0=a_1 \geq 0$, then $s=a_0$ and $d=0$.
\item If $a_0<0\le a_1$, then $s=a_1$ and $d=a_1-a_0$.
\item If $a_1<0\le a_0$, then $s=a_0$ and $d=a_0-a_1$.
\item If $a_0,a_1<0$. Here, we let $t=|a_0-a_1|$ and $m=\max(a_0,a_1)$. Set
$$
r=\min \{\ell  :  m+\ell t+T_\ell \geq 0\},
$$
so
$$
r=\left\lceil   \frac{-1-2t+\sqrt{(2t+1)^2-8m}}{2}   \right\rceil.
$$
Finally, we let $s=m+rt+T_r$ and $d=t+r$.
\end{enumerate}

Note that $s=d$ occurs only in case (3) with $s=d=0$.

\begin{thm}
\label{Main3}
Suppose that $(a_0,a_1)$ has label $(s,d)$ with $d>0$.
If $s+k d \ge T_d$, then $(a_0,a_1)$ is in 
the cycle $\mathcal{C}_{s,d,k}$ and it occurs when
$$
\lambda\in
Z_k := \left[-2+\frac {1}{s+(k+1)d}, -2+\frac{1}{s+kd} \right[.
$$
Putting $K=\lceil (T_d-s)/d\rceil$, we have
$$
\bigcup_{k \geq K}^\infty Z_k =\left]-2,  -2+ \frac{1}{s+Kd}\right[,
$$
i.e., the orbit of $(a_0,a_1)$ is periodic when $\lambda\in\ 
]-2, -2+1/(s+Kd)[$.
Suppose that $d=0$. If $s=0$, then the cycle $(0)^{\infty}$ occurs on all of $]-2,2[$. If $s>0$, then the cycle $(s)$ occurs in the interval $]-2,-2+1/s[$.
\end{thm}

\begin{proof}
For $d>0$, we can confirm that by choice of $s, d, K$, 
the initial values of $(a_0,a_1)$ occurs 
in the cycle $\mathcal{C}_{s,d,k}:=(AB\bar{B}\bar{A}C\bar{C})^{\infty}$.
In particular, in the first case, the occurrence is in $A$, in the second case, the occurrence is in $\bar{A}$. Recalling the overlapping structure of
$\mathcal{C}_{s,d,k}$ in the proof of Lemma \ref{Hex}, in the fourth case, 
the occurrence is on the boundary of $\bar{C}$ and $A$, in the fifth case, the occurrence is on the boundary of $\bar{A}$ and $C$. 
For the sixth case, by definition of $s$, we have
$$
m+(r-1)t+ T_{r-1}<0\le m+rt+ T_{r}.
$$
Thus $d-s=(t+r)-(m+rt+T_r)=-(m+(r-1)t+ T_{r-1})>0$, which guarantee the condition $0\le s<d$. The $r$-th element of $C$ is
$$
s-(T_d-T_{d-r})=m+rt+ T_{r}-(T_d-T_{d-r})=m=\max\{a_0,a_1\}<0
$$
and the difference between $(r+1)$-th element and $r$-th element
is
$$
T_{d-r}-T_{d-r-1}=d-r=t=|a_0-a_1|.
$$
Therefore $(a_0,a_1)$ appears in $C\overline{C}$.
Thus the statement follows for $k\ge K$ by Lemma \ref{Hex}.

Now, we suppose that $d=0$. If $s=0$, then $0 \leq a_2+\lambda 0+0<1$, so $a_2=0$. Next, suppose that $d>0$.
Here, $0 \leq a_2+\lambda s+s <1$. Also, assume that $-2<\lambda < -2+\frac{1}{s}$, so $0 \leq a_2-s \implies a_2 \geq s$. Next,
$$
a_2+\lambda a_1+a_0 < a_2-s+1 \leq 1 \implies a_2 \leq s.
$$
So, $a_2=s$ and $]-2,-2+1/s[$ contains only the cycle $(s)^{\infty}$.
\end{proof}

\section{Proof of Theorem \ref{Small}}

Given $(a_0,a_1)\in \Z^2$, we proved that $]-2,2[$ is partitioned into at most countably many rational intervals as in (\ref{Part}).
Further Theorem \ref{Main3} gives a complete picture around $-2$. Indeed 
it is covered by infinitely many concrete intervals. Therefore even if
Theorem \ref{Main} is a conditional result, 
if we could give a concrete partition of $[-2+\eps,2[$,
we would obtain a full description of (\ref{Part}). 
Given $(a_0,a_1)\neq (0,0)$, our algorithm is described in this way.
\begin{enumerate}
\item Compute $s,d, K$ from $(a_0,a_1)$.
\item $L\leftarrow \{1/(2s+2Kd)\}$. $X\leftarrow [-2+1/(s+Kd),2[$
\item Compute the cycle $c(\ell)$ which corresponds to $\ell\in L$ 
and give the non-empty
rational interval $I_{c(\ell)}$ which correspond to $c(\ell)$.
\item $X\leftarrow X\setminus \bigcup_{\ell\in L} I_{c(\ell)}$.
\item If $X=\emptyset$, then go to End.
\item Noting that $X$ is a union of rational intervals, 
let $L$ be the set of all mid points of the intervals of $X$. Go to (3).
\item End.
\end{enumerate}

According to Theorem \ref{Main}, this algorithm is expected to terminate. If 
it terminates for an initial vector $(a_0,a_1)\in \Z^2$, 
Conjecture \ref{SRS} is valid for $(a_0,a_1)$. 
We confirmed all values of $\max\{|a_0|,|a_1|\}\le 10$. 
End points of the 
partition of $[-2+1/(s+Kd),2)$ for $\max\{|a_0|,|a_1|\}\le 2$ are listed below.
The arrow over a rational number shows that it belongs to a proper 
interval containing this direction, and a rational number without an 
arrow is a singleton. For 
$(a_0,a_1)=(-1,1)$, e.g., the interval 
$\left]-3/2,-4/3\right[$ gives the cycle $(-1, 1, 3, 4, 3, 1, -1, -2)^{\infty}$ and 
$\left[-4/3, -1\right[$ 
gives $(-1, 1, 3, 3, 1, -1, -2)^{\infty}$, but a singleton
$[-1]$ corresponds to $(-1, 1, 2, 1, -1, -2)^{\infty}$
and $]-1,-1/2[$ gives $(-1, 1, 2, 1, -1)^{\infty}$. A 
closed interval $[-2/3,-1/2]$ appears when $(a_0,a_1)=(-2,-2)$ and it 
corresponds to  $(-2, -2, 1, 3, 1)^{\infty}$.


\begin{small}

\begin{align*}
&(a_0,a_1)=(-1,-1),\\
& -1,-\frac12,0,\frac12,1,\frac43,\frac32,\frac85,\frac53,\frac74,\frac95,2\\
&(a_0,a_1)=(-1,1),(1,-1)\\
&\stackrel{\longrightarrow}{-\frac53},-\frac32,\stackrel{\longrightarrow}{-\frac43},-1,-\frac12,0,1,\frac32,2\\
&(a_0,a_1)=(1,0),(-1,0),(0,1),(0,-1)\\
& -1,0,1,\frac32,2\\
&(a_0,a_1)=(1,1),\\
& -1,0,1,\stackrel{\longrightarrow}{\frac43},
\frac32,\frac53,2\\
&(a_0,a_1)=(-2,-2)\\
&\stackrel{\longrightarrow}{-\frac{5}{3}},-\frac{3}{2},\stackrel{\longrightarrow}{-\frac{4}{3}},-1,\stackrel{\longrightarrow}{-\frac{2}{3}},\stackrel{\longleftarrow}{-\frac{1}{2}},-\frac{1}{3},-\frac{1}{4},0,\frac{1}{4},\frac{1}{3},\stackrel{\longleftarrow}{\frac{1}{2}},\frac{2}{3},\frac{3}{4},1,\stackrel{\longrightarrow}{\frac{6}{5}},\frac{5}{4},\stackrel{\longrightarrow}{\frac{4}{3}},\frac{7}{5},\stackrel{\longleftarrow}{\frac{3}{2}},\stackrel{\longleftarrow}{\frac{11}{7}},\\
&\frac{5}{3},\stackrel{\longleftarrow}{\frac{12}{7}},\frac{7}{4},\stackrel{\longrightarrow}{\frac{16}{9}},\frac{11}{6},\stackrel{\longrightarrow}{\frac{24}{13}},\stackrel{\longleftarrow}{\frac{13}{7}},\frac{15}{8},\stackrel{\longrightarrow}{\frac{32}{17}},\stackrel{\longleftarrow}{\frac{17}{9}},\stackrel{\longrightarrow}{\frac{19}{10}},\stackrel{\longleftarrow}{\frac{21}{11}},\stackrel{\longrightarrow}{\frac{23}{12}},\stackrel{\longleftarrow}{\frac{27}{14}},\frac{29}{15},\frac{31}{16},\frac{33}{17},\frac{35}{18},2\\
&(a_0,a_1)=(-1,-2),(-2,-1)\\
&\stackrel{\longrightarrow}{-\frac{5}{3}},-\frac{3}{2},\stackrel{\longrightarrow}{-\frac{4}{3}},-1,-\frac{2}{3},-\frac{1}{2},-\frac{1}{3},0,\frac{1}{3},\frac{1}{2},\frac{2}{3},1,\frac{5}{4},\frac{4}{3},\frac{3}{2},\frac{8}{5},\frac{5}{3},\\&\stackrel{\longleftarrow}{\frac{12}{7}},\frac{19}{11},\frac{7}{4},\stackrel{\longrightarrow}{\frac{9}{5}},\stackrel{\longleftarrow}{\frac{11}{6}},\frac{13}{7},\stackrel{\longrightarrow}{\frac{15}{8}},\frac{17}{9},\frac{19}{10},\frac{21}{11},2\\
&(a_0,a_1)=(0,\pm2),(\pm2,0)\\
&\stackrel{\longrightarrow}{-\frac{7}{4}},\stackrel{\longleftarrow}{-\frac{5}{3}},\stackrel{\longrightarrow}{-\frac{8}{5}},-\frac{3}{2},\stackrel{\longrightarrow}{-\frac{4}{3}},-1,-\frac{2}{3},-\frac{1}{2},0,\frac{1}{2},1,\stackrel{\longrightarrow}{\frac{4}{3}},\frac{3}{2},\frac{8}{5},\frac{5}{3},\frac{7}{4},\frac{9}{5},\frac{11}{6},2\\
\end{align*}

\begin{align*}
&(a_0,a_1)=(1,-2),(-2,1)\\
&\stackrel{\longrightarrow}{-\frac{13}{7}},-\frac{9}{5},\stackrel{\longrightarrow}{-\frac{16}{9}},-\frac{7}{4},\stackrel{\longrightarrow}{-\frac{12}{7}},\stackrel{\longrightarrow}{-\frac{5}{3}},-\frac{3}{2},\stackrel{\longrightarrow}{-\frac{7}{5}},-\frac{4}{3},\stackrel{\longrightarrow}{-\frac{5}{4}},-1,\stackrel{\longrightarrow}{-\frac{2}{3}},\stackrel{\longleftarrow}{-\frac{1}{2}},-\frac{1}{3},0,\frac{1}{2},1,\stackrel{\longrightarrow}{\frac{4}{3}},\frac{3}{2},\frac{5}{3},2\\
&(a_0,a_1)=(2,-2),(-2,2)\\
&\stackrel{\longrightarrow}{-\frac{19}{10}},-\frac{15}{8},\stackrel{\longrightarrow}{-\frac{28}{15}},\stackrel{\longleftarrow}{-\frac{13}{7}},\stackrel{\longrightarrow}{-\frac{24}{13}},\stackrel{\longrightarrow}{-\frac{11}{6}},\stackrel{\longrightarrow}{-\frac{20}{11}},\stackrel{\longleftarrow}{-\frac{9}{5}},\stackrel{\longrightarrow}{-\frac{16}{9}},\stackrel{\longrightarrow}{-\frac{17}{10}},-\frac{5}{3},\stackrel{\longrightarrow}{-\frac{13}{8}},\stackrel{\longleftarrow}{-\frac{8}{5}},-\frac{11}{7},\\
&\stackrel{\longrightarrow}{-\frac{14}{9}},-\frac{3}{2},\stackrel{\longrightarrow}{-\frac{10}{7}},\stackrel{\longrightarrow}{-\frac{7}{5}},-\frac{4}{3},-\frac{5}{4},\stackrel{\longrightarrow}{-\frac{6}{5}},-1,\stackrel{\longrightarrow}{-\frac{3}{4}},-\frac{2}{3},-\frac{1}{2},-\frac{1}{3},-\frac{1}{4},0,\frac{1}{3},\frac{1}{2},1,\stackrel{\longrightarrow}{\frac{4}{3}},\frac{3}{2},2\\
&(a_0,a_1)=(2,-1),(-1,2)\\
&\stackrel{\longrightarrow}{-\frac{15}{8}},\stackrel{\longrightarrow}{-\frac{20}{11}},\stackrel{\longrightarrow}{-\frac{9}{5}},-\frac{7}{4},\stackrel{\longrightarrow}{-\frac{12}{7}},\stackrel{\longleftarrow}{-\frac{5}{3}},
\stackrel{\longrightarrow}{-\frac{8}{5}},\stackrel{\longrightarrow}{-\frac{3}{2}},\stackrel{\longleftarrow}{-\frac{4}{3}},\stackrel{\longrightarrow}{-\frac{5}{4}},-1,-\frac{2}{3},-\frac{1}{2},-\frac{1}{3},0,\frac{1}{2},1,\frac{3}{2},2\\
&(a_0,a_1)=(2,1),(1,2)\\
&-1,-\frac{1}{2},0,\frac{1}{2},1,\frac{5}{4},\frac{4}{3},\frac{3}{2},\stackrel{\longrightarrow}{\frac{8}{5}},\stackrel{\longleftarrow}{\frac{5}{3}},\stackrel{\longleftarrow}{\frac{12}{7}},\frac{19}{11},\frac{7}{4},\stackrel{\longrightarrow}{\frac{9}{5}},\frac{11}{6},\frac{13}{7},\frac{15}{8},2\\
&(a_0,a_1)=(2,2)\\
&\stackrel{\longrightarrow}{-\frac{3}{2}},-1,-\frac{2}{3},-\frac{1}{2},-\frac{1}{3},0,\frac{1}{3},\frac{1}{2},\frac{2}{3},1,\stackrel{\longrightarrow}{\frac{6}{5}},\frac{5}{4},\frac{4}{3},\frac{3}{2},\frac{8}{5},
\frac{13}{8},\frac{5}{3},\frac{12}{7},
\stackrel{\longleftarrow}{\frac{9}{5}},\stackrel{\longrightarrow}{\frac{20}{11}},\stackrel{\longleftarrow}{\frac{13}{7}},\stackrel{\longrightarrow}{\frac{15}{8}},\\
&\stackrel{\longleftarrow}{\frac{19}{10}},\frac{21}{11},\frac{23}{12},\frac{25}{13},\frac{27}{14},2
\end{align*}
\end{small}

The maximum number of intervals in a partition for a fixed $m$ is listed below.
The cardinality of the partitions may be of quadratic growth.

$$
\begin{tabular}{l*{6}{c}r}
$m$ & Max at $(a_0,a_1)$              & Cardinality of Partition & Cardinality of Singletons   \\
\hline
$1$ & $(-1,-1)$ & $22$ & $11$ \\
$2$ & $(-2,-2)$ & $59$ & $20$ \\
$3$ & $(-3,-3)$ & $135$ & $52$ \\
$4$ & $(-3,-4)$ & $222$ & $92$ \\
$5$ & $(-4,-5)$ & $369$ & $147$ \\
$6$ & $(-5,-6)$ & $520$ & $209$ \\
$7$ & $(-6,-7)$ & $674$ & $260$ \\
$8$ & $(-7,-8)$ & $956$ & $388$ \\
$9$ & $(-8,-9)$ & $1143$ & $442$ \\
$10$ & $(-9,-10)$ & $1409$ & $554$ \\
\end{tabular}
$$

The maximum cycle length as well as
the average cycle length are listed below. Average lengths may be of 
linear growth. Thus the number of applications of the map
would be $O(m)\times O(m^2)=O(m^3)$.

$$
\begin{tabular}{ccccc}
$m$ & Max at $(a_0,a_1)$ & Interval & Maximum length & Average length   \\
\hline
$1$ & $(-1,-1)$ & $[8/5]$ & $38$ & $9.8172$ \\
$2$ & $(2,2)$ & $[25/13]$ & $123$ & $21.5067$ \\
$3$ & $(3,3)$ & $[11/6]$ & $253$ & $33.9038$ \\
$4$ & $(4,4)$ & $]15/8,62/33[$ & $363$ & $50.7101$ \\
$5$ & $(-5,-5)$ & $[137/70]$& $623$ & $67.6183$ \\
$6$ & $(-5,-6)$ & $]139/71,47/24[$ & $952$ & $85.2249$ \\
$7$ & $(-7,-7)$ & $[75/38]$ & $1065$ & $99.7307$ \\
$8$ & $(-7,-8)$ & $[391/196]$ & $1040$ & $120.428$ \\
$9$ & $(-8,-9)$ & $]145/73,292/147[$ & $1365$ & $139.302$ \\
$10$ & $(10,9)$ &$[78/41]$ & $1802$ & $158.299$ \\
\end{tabular}
$$
From this numerical observation, we guess that our algorithm requires $O(m^3)$ operations to determine all the partitions for all $(a_0,a_1)$ with $m=\max(|a_0|,|a_1|)$. We expect determining all partitions for $(a_0,a_1)$ with $\max(|a_0|,|a_1|) \leq m$ to be $O(m^4)$. This conjecture is numerically reconfirmed 
from the file sizes that our program created for a fixed value $m$. 

\end{document}